\documentclass[a4paper]{amsart}

\usepackage{amsthm, amsfonts, amssymb, amsmath, latexsym, enumerate, times}
\usepackage[latin1]{inputenc}
\usepackage{mathrsfs}
\usepackage{mathtools}
\usepackage{letltxmacro}
\LetLtxMacro\orgvdots\vdots
\LetLtxMacro\orgddots\dots

\makeatletter
\DeclareRobustCommand\vdots{%
	\mathpalette\@vdots{}%
}
\newcommand*{\@vdots}[2]{%
	\sbox0{$#1\cdotp\cdotp\cdotp\m@th$}%
	\sbox2{$#1.\m@th$}%
	\vbox{%
		\dimen@=\wd0 %
		\advance\dimen@ -3\ht2 %
		\kern.5\dimen@
		\dimen@=\wd2 %
		\advance\dimen@ -\ht2 %
		\dimen2=\wd0 %
		\advance\dimen2 -\dimen@
		\vbox to \dimen2{%
			\offinterlineskip
			\copy2 \vfill\copy2 \vfill\copy2 %
		}%
	}%
}
\DeclareRobustCommand\ddots{%
	\mathinner{%
		\mathpalette\@ddots{}%
		\mkern\thinmuskip
	}%
}
\newcommand*{\@ddots}[2]{%
	\sbox0{$#1\cdotp\cdotp\cdotp\m@th$}%
	\sbox2{$#1.\m@th$}%
	\vbox{%
		\dimen@=\wd0 %
		\advance\dimen@ -3\ht2 %
		\kern.5\dimen@
		\dimen@=\wd2 %
		\advance\dimen@ -\ht2 %
		\dimen2=\wd0 %
		\advance\dimen2 -\dimen@
		\vbox to \dimen2{%
			\offinterlineskip
			\hbox{$#1\mathpunct{.}\m@th$}%
			\vfill
			\hbox{$#1\mathpunct{\kern\wd2}\mathpunct{.}\m@th$}%
			\vfill
			\hbox{$#1\mathpunct{\kern\wd2}\mathpunct{\kern\wd2}\mathpunct{.}\m@th$}%
		}%
	}%
}
\makeatother

\usepackage{array}
\usepackage{comment}
\usepackage{paralist}
\usepackage{xcolor}

\newtheorem{theorem}{Theorem}
\newtheorem{lemma}[theorem]{Lemma}

\newtheorem{corollary}[theorem]{Corollary}

\newtheorem{question}[theorem]{Question}

\theoremstyle{definition}

\newtheorem{example}[theorem]{Example}

\newcommand{\bbQ}{{\mathbb Q}}
\newcommand{\bbC}{{\mathbb C}}

\def\le{\leqslant}
\def\ge{\geqslant}

\usepackage{eurosym}

\begin{document}
 
\title{A remark on generalized abundance for surfaces}
 
\author{Claudio Fontanari}

\address{Claudio Fontanari, Dipartimento di Matematica, Universit\`a degli Studi di Trento, Via Sommarive 14, 38123 Povo, Trento}
\email{claudio.fontanari@unitn.it}

\subjclass{14E30} 
 
\keywords{Abundance, Generalized Abundance, Semiampleness, Numerical Semiampleness}
 
\centerline{}
\begin{abstract} 
Let $(X, \Delta)$ be a projective klt pair of dimension $2$ and let $L$ be a nef Cartier divisor on $X$ such that $K_X + \Delta + L$ is nef. 
As a complement to the Generalized Abundance Conjecture by Lazi\'c and Peternell, we prove that 
if $K_X + \Delta$ and $L$ are not proportional modulo numerical equivalence, then $K_X + \Delta + L$ is semiample. 
An example due to Lazi\'c shows that this is no longer true in any dimension $n \ge 3$. 
\end{abstract}

\maketitle 

\section{Introduction}

The Generalized Abundance Conjecture by Lazi\'c and Peternell (see \cite{LP}, p. 354) is indeed a theorem in dimension $2$ 
(see \cite{LP}, Corollary C. on p. 356): 

\begin{theorem} 
Let $(X, \Delta)$ be a projective klt pair of dimension $2$ such that $K_X + \Delta$ is pseudoeffective and let $L$ be a nef Cartier divisor on $X$. If $K_X + \Delta + L$ is nef then there exists a semiample $\bbQ$-divisor $M$ on $X$ such that 
$K_X + \Delta + L$ is numerically equivalent to $M$. 
\end{theorem}

The assumption that $K_X + \Delta$ pseudoeffective turns out to be necessary (see for instance \cite{LP}, Example 6.2). On the 
other hand, at least in dimension $2$, it is possible to characterize the failure of numerical abundance when $K_X + \Delta$ is not pseudoeffective. The following statement is \cite{HL}, Theorem 3.13: 

\begin{theorem} 
Let $(X, \Delta)$ be a projective klt pair of dimension $2$ and let $L$ be a nef Cartier divisor on $X$ such that $K_X + \Delta + L$ is nef. Then either $K_X + \Delta + L$ is numerically semiample or $K_X + \Delta$ is numerically equivalent to $-tL$ with $0 \le t \le 1$. 
\end{theorem}

We point out that if $t=0$ then we fall in the first case. Indeed, the Semiampleness Conjecture holds on surfaces (see \cite{LP}, Theorem 8.2): 
if $K_X + \Delta$ is numerically equivalent to $0$ then $L$ is numerically semiample. It is therefore tempting to ask the following question in higher dimension: 
 
\begin{question}
Let $(X, \Delta)$ be a projective klt pair of dimension $n \ge 3$ and let $L$ be a nef Cartier divisor on $X$ such that $K_X + \Delta + L$ is nef. Is it true that either $K_X + \Delta + L$ is numerically semiample or $L$ is numerically equivalent to $-m(K_X + \Delta)$ with $m > 0$?  
\end{question}

Even though we are not aware of any counterexamples, there seems to be no reason to expect an affirmative answer. 

As shown already in dimension $1$ by the example of a non-torsion numerically trivial divisor on an elliptic curve (see \cite{BH}, p. 212), numerical semiampleness cannot be replaced by semiampleness. We notice however that, at least in dimension $2$, semiampleness holds under an easily stated explicit assumption. We formulate this remark as follows:

\begin{theorem}\label{main}
Let $(X, \Delta)$ be a projective klt pair of dimension $2$ and let $L$ be a nef Cartier divisor on $X$ such that $K_X + \Delta + L$ is nef. 
If $K_X + \Delta$ and $L$ are not proportional modulo numerical equivalence, then $K_X + \Delta + L$ is semiample. 
\end{theorem}

The above result complements but does not imply Generalized Abundance, in particular its statement is empty in the two crucial cases 
$L = K_X + \Delta$ (Abundance Conjecture) and $K_X + \Delta$ numerically trivial (Semiampleness Conjecture on Calabi-Yau pairs).
Once again, it is legitimate to wonder about the higher dimensional case.  We are going to present our proof in a general setting, but in arbitrary dimension we only obtain a pale shadow of the two dimensional case (see Corollary \ref{proportional}). The following example, kindly provided to us by Vladimir Lazi\'c, shows that the statement of Theorem \ref{main} does not extend to any dimension $n \ge 3$: 

\begin{example} \textbf{(Lazi\'c)} Let $X$ be a smooth variety with $\mathrm{Pic}^0(X)=0$ and Picard number at least $2$. Take an ample divisor $A$ on $X$ not proportional to $K_X$ and such that $K_X+A$ is ample. Let $E$ be an elliptic curve and take a degree zero non-torsion divisor $P$ on $E$. Consider $Y=X \times E$, let $A_Y$ be the pullback of $A$ to $Y$ via the first projection and let $P_Y$ be the pullback of $P$ via the second projection. Then $K_Y$ and $A_Y+P_Y$ are not proportional modulo numerical equivalence, but $K_Y+A_Y+P_Y$ is not semiample. Indeed, assume by contradiction that $K_Y+A_Y+P_Y$ is semiample and consider the induced Iitaka fibration $f: Y\to Z$. Then $f$ and the first projection $Y \to X$ contract the same curves, hence $X$ and $Z$ are isomorphic by the rigidity lemma. From the factorization $f: Y\to X \to Z$ it follows that $P_Y$ is the pullback (up to $\bbQ$-linear equivalence) of a divisor $P_X$ from $X$, since $K_Y+A_Y$ is the pullback of a divisor from $X$ and $K_Y+A_Y+P_Y$ is the pullback of a divisor from Z.
Then $P_X$ is numerically trivial on $X$, hence torsion by the assumption $\mathrm{Pic}^0(X)=0$. But this would imply that $P_Y$ is torsion, hence $P$ is torsion, a contradiction.
\end{example}

We work over the complex field $\bbC$.

\medskip

{\bf Acknowledgements:} 
The author is grateful to Vladimir Lazi\'c and Thomas Peternell for their helpful comments. 
This research project was partially supported by GNSAGA of INdAM and by PRIN 2017 ``Moduli Theory and Birational Classification''. 

\section{The proof}

Our first Lemma generalizes \cite{S}, Lemma 1.3.

\begin{lemma} \label{big}
Let $(X, \Delta)$ be a projective klt pair of dimension $n$ and let $H$ be a nef and big Cartier divisor on $X$. 
If $L$ is a nef Cartier divisor on $X$ such that $K_X + \Delta + L$ is nef and $K_X + \Delta + 2 L$ has 
numerical dimension $\nu(K_X + \Delta + 2 L) < k \le n$, then we have 
$$
H^{n-k}L^k = H^{n-k} L^{k-1}(K_X + \Delta) = \ldots = H^{n-k} (K_X + \Delta)^k = 0.
$$

\end{lemma}

\proof Since both $K_X + \Delta + L$ and $L$ are nef we have 
$$
0 \le H^{n-k}(K_X + \Delta + 2L)^k = \sum_{m=0}^k \binom{k}{m}H^{n-k}(K_X + \Delta + L)^m L^{k-m}
$$
with $H^{n-k}(K_X + \Delta + L)^m L^{n-m} \ge 0$ for every $m$. 

If $H^{n-k}(K_X + \Delta + 2L)^k = 0$ then $H^{n-k}(K_X + \Delta + L)^m L^{k-m} = 0$ for every $m$ and by induction it follows that 
$H^{n-k}L^k = H^{n-k}L^{k-1}(K_X + \Delta) = \ldots = H^{n-k}(K_X + \Delta)^k = 0$.

\qed

\begin{corollary} \label{semiample}
Let $(X, \Delta)$ be a projective klt pair of dimension $n$ and let $L$ be a nef Cartier divisor on $X$. If $K_X + \Delta + L$ is nef but not semiample, 
then we have 
$$
L^n = L^{n-1}(K_X + \Delta) = \ldots = (K_X + \Delta)^n = 0.
$$
\end{corollary}

\proof 

If $(K_X + \Delta + 2L)^n > 0$ then $K_X + \Delta + 2L = 2(K_X + \Delta + L) - (K_X + \Delta)$ is nef and big, hence 
$K_X + \Delta + L$ would be semiample by the logarithmic base-point-free theorem. Since $K_X + \Delta + L$ is not semiample
we deduce that $(K_X + \Delta + 2L)^n = 0$ and $\nu(K_X + \Delta + 2 L) < n$. Now the claim follows from Lemma \ref{big} with $k=n$. 

\qed

Our next Lemma generalizes to arbitrary dimension the \textit{Easy fact} stated for surfaces in \cite{CCP}, pp. 576-577 (see also \cite{HL}, 
Lemma 3.2, where the assumption $A^2=B^2=0$ is missing and the assumption $A$, $B$ nef is added). 

\begin{lemma}\label{index}
Let $X$ be a normal projective variety of dimension $n$ and let $H$ be a nef and big Cartier divisor on $X$. 
If $A$ and $B$ are two $\bbQ$-Cartier divisors on $X$ such that $H^{n-2}A^2 = H^{n-2}B^2 = H^{n-2}AB = 0$, 
then $A$ and $B$ are proportional modulo numerical equivalence.  
\end{lemma}

\proof By replacing $X$ with a birational resolution of singularities and $A$ and $B$ by their pullbacks 
we may assume that $X$ is smooth. We may also assume that $H^{n-1}A$ and $H^{n-1}B$  
are proportional by a rational factor $m$, so that $H^{n-1}(A-mB)=0$. Now we apply the Hodge index theorem for divisors 
(see \cite{T}, \S 1, and \cite{L}, Theorem 1) to $E = A-mB$: if $H^{n-1}E=0$ then $H^{n-2}E^2 \le 0$
and equality holds if and only if $H^{n-2}E$ is homologically equivalent to zero. By assumption we 
have $H^{n-2}E^2 = H^{n-2}(A-mB)^2 =  H^{n-2}A^2 + m^2 H^{n-2}B^2 - 2m H^{n-2}AB = 0$, hence 
$H^{n-2}E = H^{n-2}(A-mB)$ is homologically equivalent to zero. By the hard Lefschetz theorem
(see for instance \cite{M}, Theorem 4.6), the Lefschetz operator $H^{n-2}: H^2(X, \bbQ) \to H^{2n-2}(X, \bbQ)$ 
is injective, therefore $A-mB$ is homologically (in particular, numerically) equivalent to zero.

\qed

\begin{corollary}\label{proportional}
Let $(X, \Delta)$ be a projective klt pair of dimension $n$. 
If $L$ is a nef divisor on $X$ such that $K_X + \Delta + L$ is nef and $K_X + \Delta + 2 L$ has 
numerical dimension $\nu(K_X + \Delta + 2 L) < 2$, then $K_X + \Delta$ 
and $L$ are proportional modulo numerical equivalence.
\end{corollary}

\proof Let $H$ be a nef and big Cartier divisor on $X$. By Lemma \ref{big} with $k=2$ we have 
$H^{n-2}L^2 = H^{n-2} L(K_X + \Delta) = H^{n-2} (K_X + \Delta)^2 = 0$.
Now the claim follows from Lemma \ref{index}.

\qed

\noindent \textit{Proof of Theorem \ref{main}.}
We argue by contradiction. If $K_X + \Delta + L$ is not semiample, then by Corollary \ref{semiample} 
we have $L^2 = L(K_X + \Delta) = (K_X + \Delta)^2 = 0$. 
Now Lemma \ref{index} yields the sought-for contradiction. 

\qed

\end{document}